\documentclass[11pt]{article}
 \usepackage{amsmath,amssymb}
 \usepackage{epsfig}
 \usepackage{graphicx}
 \usepackage{pict2e}
\usepackage{verbatim}
\usepackage{multirow} 
\usepackage{bm}
 \newtheorem{Lemma}{Lemma}
 \newtheorem{Proposition}[Lemma]{Proposition}

 \newtheorem{Conjecture}[Lemma]{Conjecture}

 \newtheorem{Corollary}[Lemma]{Corollary}

 \newcommand{\FF}{\mbox{${\mathcal F}$}} 
  \newcommand{\GG}{\mbox{${\mathcal G}$}} 
  
 \renewcommand{\SS}{\mbox{${\mathcal S}$}}
 
 \newcommand{\sfrac}[2]{{\textstyle\frac{#1}{#2}}}

  \newcommand{\Ints}{{\mathbb{Z}}}
  \newcommand{\bw}{{\mathbf w}}
    \newcommand{\bV}{{\mathbf V}}
  \newcommand{\bE}{{\mathbf E}}

\def\Ex{\mathbb{E}}
\def\Pr{\mathbb{P}}
\def\ind{{\rm 1\hspace{-0.90ex}1}}

\newcommand{\var}{\mathrm{var}}

\begin{document}

\title{The Incipient Giant Component in Bond  Percolation on General Finite Weighted Graphs}
 \author{David J. Aldous
 \thanks{Department of Statistics,
 367 Evans Hall \#\  3860,
 U.C. Berkeley CA 94720;  aldous@stat.berkeley.edu;
  www.stat.berkeley.edu/users/aldous.  Aldous's research supported by
 N.S.F Grant DMS-1504802. }}

 \maketitle
 
  \begin{abstract}
On a large finite connected graph let edges $e$ become ``open" at independent random Exponential times of arbitrary rates $w_e$. 
Under minimal assumptions, the time at which a giant component starts to emerge is weakly concentrated around its mean.
  \end{abstract}
 \vspace{0.1in}


\section{Introduction}
\label{sec:intro}
 Take a finite connected graph $(\bV,\bE)$ with edge-weights $\bw = (w_e)$, where $w_e > 0 \ \forall e \in \bE$.
 To the edges $e \in \bE$ attach independent Exponential(rate $w_e$) random variables $\xi_e$.
 In the language of percolation theory, say that edge $e$ becomes {\em open} at time $\xi_e$. 
 The set of open edges at time $t$ determines a random partition of $\bV$ into connected components; 
 write  $C(t)$ for the largest number of vertices in any such connected component.  
 Now consider a sequence $(\bV_n,\bE_n)$ of such weighted  graphs, where both the graph topologies and the edge-weights 
 are arbitrary subject only to the conditions that $|\bV_n| \to \infty$ and that for some $0 < t_1 < t_2 < \infty$ 
 \begin{equation}
 \lim_n  \Ex C_n(t_1)/|\bV_n| = 0; \quad  \bar{c} := \liminf_n  \Ex C_n(t_2)/|\bV_n| > 0 . 
 \label{hypo}
 \end{equation}
 In the language of random graph theory, this condition says that a {\em giant component} emerges (with non-vanishing probability) sometime between $t_1$ and $t_2$.
 Proposition \ref{P1} asserts, informally,  that the ``incipient" time at which the giant component starts to emerge is deterministic to first order.  
 \begin{Proposition}
 \label{P1}
 Given a sequence of graphs satisfying (\ref{hypo}), there exists a deterministic sequence $\tau_n \in [t_1,t_2]$ such that, for every 
 sequence  $\omega_n \uparrow \infty$ sufficiently slowly, the random times
 \[ T_n := \inf \{t: C_n(t) \ge |\bV_n|/\omega_n\} \]
 satisfy
 \[ T_n - \tau_n \to_p 0 . \]
 \end{Proposition}
 
In the special cases of the complete graph and the 2-dimensional discrete torus (with constant edge-weights) 
we are essentially dealing with component sizes in the classical
 {E}rd{\H o}s-{R}\'enyi  and  the bond percolation on $\Ints^2$ processes, 
 for which much stronger results are known about the
``scaling window" of time over which the giant component emerges 
\cite{boll-survey,borgs}. 
Such stronger results have been generalized (again with constant edge-weights) in several directions, for instance to random subgraphs of certain transitive finite graphs 
\cite{CJH-1,CJH-2} or to random subgraphs of graphs under assumptions that force the critical subgraphs to be ``tree-like" 
as in the  {E}rd{\H o}s-{R}\'enyi case \cite{chung-horn}. 
 Proposition \ref{P1} gives a comparatively  weak concentration property which one would expect to hold in every ``natural" example,
 but it is perhaps remarkable that it holds in the generality stated.  
 As a comparison, in the classical cases one also has the same ``weak concentration" property for the random times 
  \[ T^*_n(s)  := \inf \{t: C_n(t) \ge  s  |\bV_n| \} \]
  for fixed $0 < s < 1$. 
  One would expect this to extend to other ``natural" examples, but a simple example outlined in section \ref{sec:example} shows it does not hold 
  in the generality of  Proposition \ref{P1}, even if we assume the weighted graph to be vertex-transitive.  
  Some conjectures concerning the post-incipient regime are given in section \ref{sec:conjs}.

The proof is based on a simple general variance bound described in section \ref{sec:gvb}.
A technically more complicated application of that bound to first passage percolation on general weighted graphs, 
plus other simple applications, can be found in \cite{me-FPPG}.
``Big picture" discussions of various random processes over finite edge-weighted graphs can be found in \cite{aldous-FMIE}
 and \cite{aldous-lisha}.

 \section{Proof of Proposition \ref{P1}}

We divide the proof into three steps.
 
 \subsection{Step 1: The general variance bound for increasing set-valued processes.} 
 \label{sec:gvb}
 Fix a weighted graph $(\bV,\bE,\bw)$ and $1< \omega < | \bV|$.
 The process of open edges in our bond percolation process is a continuous-time Markov chain, $Z_t$ say, whose state space is the set 
 of subsets $S\ \subseteq \bE$ and whose transition rates are
 \[ S \to S \cup \{e\}: \quad \mbox{ rate } w_e, \quad (e \not\in S) . \]
   We seek to study the distribution of the stopping time
 \begin{equation}
 T =  \inf \{t: C(t) \ge |\bV|/\omega\} 
 \label{T-def}
 \end{equation}
 when this chain starts in state $\emptyset$.  
 It makes sense to also consider this mean hitting time started from an arbitrary subset $S$ of open edges, that is
 \[ h(S) := \Ex_S T \]
 which clearly has the property
 \begin{equation}
h(S^\prime) \le h(S) \mbox{ whenever $S \to S^\prime$ is a possible transition}. 
\label{def-monotone}
\end{equation}
There is a general concentration inequality for Markov chain stopping times with this property: it is an easy 
consequence of martingale identities for $\Ex T$ and $\var \ T$,
though apparently not well known. 
In general the starting state is arbitrary, but to fit our setting we take it as $\emptyset$.

\begin{Lemma}[\cite{me-FPPG} equation (9) and subsequent display]
\label{L-simple-3}
Given property (\ref{def-monotone}), for arbitrary $\delta > 0$,
\[
 \frac{\var_\emptyset \ T}{(\Ex_\emptyset T)^2} \le \delta  
+\frac{  \Ex_\emptyset \int_0^T q_\delta (Z_u)  du }{\Ex_\emptyset T}  
\]
where
\begin{equation}
 q_\delta (S) :=  \sum_{S^\prime : \ h(S) - h(S^\prime) >  \delta \Ex_\emptyset T } 
\ q(S,S^\prime) (h(S) - h(S^\prime) ) \le 1
\label{def-q}
\end{equation}
and $q(S,S^\prime) $ are the transition rates.
\end{Lemma}
We will use a general consequence, derived as Corollary \ref{Cor:1} below,
Note that $h(S) \le h(\emptyset) = \Ex_\emptyset T$ for subsets $S$ under consideration.
Now consider the first time (if ever) $W_\delta$ that the process makes some transition $S \to S^\prime$ with 
 $h(S) - h(S^\prime)$ greater than $  \delta \Ex_\emptyset T$:
\begin{equation}
W_\delta := \inf \{t:  h(Z_{t-}) - h(Z_{t}) > \delta \Ex_\emptyset T \}.
\label{def-W}
\end{equation}
Then
\[  \int_0^T q_\delta (Z_u)  du =   \int_0^{T \wedge W_\delta}  q_\delta (Z_u)  du +  \int_{T \wedge W_\delta} ^T q_\delta (Z_u)  du \]
\[ \le \int_0^{T \wedge W_\delta}  q_\delta (Z_u)  du + T \ind_{\{ W_\delta<T\}} .\]
Property (\ref{def-monotone}) implies the distribution of $T$ has the submultiplicativity property
\begin{equation}
 \Pr_\emptyset(T > t_1+t_2) \le \Pr_\emptyset(T > t_1) \ \Pr_\emptyset(T > t_2) ,\quad t_1, t_2 > 0 
 \label{submult}
 \end{equation}
and it is straightforward to show 
there is a certain ``universal" function $\gamma(u) \downarrow 0$ as $u \downarrow 0$ such that, for any submultiplicative 
$T$ and any event $A$ we have 
$\Ex [T \ind_A] \le \gamma(\Pr(A)) \ \Ex T$.
So in our setting
\[ \Ex_\emptyset [T \ind_{\{ W_\delta <T\}} ] \le \gamma( \Pr_\emptyset  (W_\delta <T)) \ \Ex_\emptyset T  . \]
Then  from Lemma \ref{L-simple-3}
\[
 \frac{\var_\emptyset \ T}{(\Ex_\emptyset T)^2} \le \delta  + \gamma( \Pr_\emptyset   (W_\delta <T)) 
+\frac{  \Ex_\emptyset \int_0^{T \wedge W_\delta}  q_\delta (Z_u)  du }{\Ex_\emptyset T}  .
\]
Now consider
\[ \widetilde{q}_\delta(S) :=  \sum_{S^\prime : \ h(S) - h(S^\prime) >  \delta \Ex_\emptyset T } 
\ q(S,S^\prime) \]
so that 
\[ q_\delta(S) \le \widetilde{q}_\delta(S) \ \Ex_\emptyset T . \]
But $\widetilde{q}_\delta(Z_u) $ is the intensity rate of $W_\delta$, that is
\[ \ind_{\{ W_\delta \le t\}} - \int_0^{W _\delta \wedge t} \widetilde{q}_\delta(Z_u) du \mbox{ is a martingale}, \]
and so from the optional sampling theorem
\[ \Ex_\emptyset \int_0^{W_\delta \wedge T} \widetilde{q}_\delta(Z_u) du = \Pr_\emptyset (W_\delta \le T) . \]
Combining this with the previous two displayed inequalities gives 
\begin{Corollary}
\label{Cor:1}
\begin{equation}
 \frac{\var_\emptyset \ T}{(\Ex_\emptyset T)^2} \le \delta  + \gamma( \Pr_\emptyset   (W_\delta \le T)) 
+ \Pr_\emptyset ( W_\delta \le T) .
\label{veT}
\end{equation}
\end{Corollary}
We remark that  the argument for inequality (\ref{veT}) does not use any structure of the bond percolation process except property (\ref{def-monotone}).

\subsection{Step 2: Bounding in terms of the growth rate of the incipient giant component.}
 We now use the structure of the bond percolation process by relating $T$ defined at (\ref{T-def}) to
 \[ T^{(2)}  =  \inf \{t: C(t) \ge 2 |\bV|/\omega\}  \ge T . \]
 Consider a possible transition $S^\prime \to S^{\prime \prime} = S^\prime \cup \{e\}$.
 We will show
 \begin{equation}
 h(S^\prime) - h(S^{\prime \prime}) \le \Ex_{ S^{\prime \prime}} (T^{(2)} - T) .
 \label{hhSS}
 \end{equation}
 There is a natural coupling $(Z^\prime_u, Z^{\prime\prime}_u, u \ge 0)$ of the processes started from 
 $S^\prime$ and from $ S^{\prime \prime} $; that is, 
 $Z^{\prime\prime}_u = Z^\prime_u \cup \{e\}$ until the Exponential($w_e$) time at which $e \in Z^\prime_u$, after which time 
  $Z^{\prime\prime}_u = Z^\prime_u$.
Write $C^\prime(u), C^{\prime\prime}(u)$ for the largest component sizes, and 
 $T^\prime, T^{\prime\prime}$ for the stopping times  (\ref{T-def}), applied to these coupled processes.
 At time $T ^{(2) \prime\prime} =  \inf \{t: C^{\prime \prime}(t) \ge 2|\bV|/\omega\} $
 the process $Z^{\prime\prime}$ contains a component of size at least $2|\bV|/\omega$, 
 and so after deleting edge $e$ there must remain a component of size at least $|\bV|/\omega$ of $Z^\prime$.
 That establishes the second inequality in
 \[ T^{\prime \prime}  \le T^\prime \le T ^{(2) \prime\prime} \]
and the first inequality is immediate. 
Now
\[  h(S^\prime) - h(S^{\prime \prime}) = \Ex T^\prime - \Ex T^{\prime \prime} 
\le \Ex  T ^{(2) \prime\prime} - \Ex T^{\prime \prime} =  \Ex_{ S^{\prime \prime}} (T^{(2)} - T) \]
establishing (\ref{hhSS}).

 Now let $V$ be the time (if any) that the bond percolation process started at $\emptyset$ makes a specified transition.
 Then (\ref{hhSS}) says that
 \begin{equation}
  h(Z_{V-}) - h(Z_V) \le \Ex_\emptyset (T^{(2)} - T \vert \FF_V) \mbox{ on } \{V \le T\} . 
  \label{hZV}
  \end{equation}
 Now fix $\delta > 0$.  For each pair $\pi = (S^\prime, S^\prime \cup \{e\})$ for which 
 $h(S^\prime) - h(S^\prime \cup \{e\}) \ge \delta \Ex_\emptyset T$ 
 there is a time $V_\pi$ as above. 
 The random variable $W_\delta$ at (\ref{def-W}) is 
 $W_\delta = \min_\pi V_\pi$, and so (\ref{hZV}) implies
  \begin{equation}
  \delta \Ex_\emptyset T \le 
  h(Z_{W_\delta -}) - h(Z_{W_\delta} ) \le \Ex_\emptyset (T^{(2)} - T \vert \FF_{W_\delta}) \mbox{ on } \{W_\delta \le T\} . 
  \label{hZW}
  \end{equation}
  This in turn implies
  \[ \Pr_\emptyset (W_\delta \le T) \le \frac{ \Ex_\emptyset (T^{(2)} - T)}{   \delta \Ex_\emptyset T} . \]
 Applying (\ref{veT}), and setting
 \[ \gamma^*(u) = \gamma(u) + u \downarrow 0 \mbox{ as } u \downarrow 0 , \]
 we find
 \[
 \frac{\var_\emptyset \ T}{(\Ex_\emptyset T)^2} \le \delta  + \gamma^*\left(  \frac{ \Ex_\emptyset (T^{(2)} - T)}{   \delta \Ex_\emptyset T}       \right) 
\]
 Because $\delta$ is arbitrary, this implies
  \begin{equation}
 \frac{\var_\emptyset \ T}{(\Ex_\emptyset T)^2} \le \Gamma \left(  \frac{ \Ex_\emptyset (T^{(2)} - T)}{ \Ex_\emptyset T}       \right) 
\label{veW}
\end{equation}
where
\[ \Gamma(x) := \inf_{\delta > 0} ( \delta + \gamma^*(x/\delta) ) \downarrow 0 \mbox{ as } x \downarrow 0 .
  \]

 \subsection{Step 3: A compactness reduction.}
 \label{sec:compact}
 The remainder of the proof uses only ``soft" arguments.
 We are given a sequence of weighted graphs satisfying (\ref{hypo}).
 To emphasize dependence on $\omega_n$ write
  \[ T_n(\omega_n) := \inf \{t: C_n(t) \ge |\bV_n|/\omega_n\} . \]
 Take $\omega_n \ge 2$ to avoid trivialities.  
  By the second condition in (\ref{hypo}) and submultiplicativity (\ref{submult}) there is an integrable $T^*$ such that
 \begin{equation}
 \mbox{
 $T^*$ stochastically dominates $T_n(\omega_n)$, for all $n, \  \omega_n$} .
 \label{dom}
 \end{equation}
 By the first condition in (\ref{hypo}) we can take $\omega_n \uparrow 0$ sufficiently
  slowly that $\Pr (T_n \le  t_1) \to 0$.  
  Looking at (\ref{veW}), we see that the proof of Proposition \ref{P1} reduces to the proof of
  \begin{quote}
  (*) for all  $\omega_n \uparrow \infty$ sufficiently slowly, 
  $\Ex ( T_n(\sfrac{1}{2} \omega_n) -  T_n(\omega_n)) \to 0$.
  \end{quote}

  Property (\ref{dom}) implies compactness with respect to weak convergence and convergence of expectations.
  By a standard compactness principle, to prove (*) it will suffice to prove that every subsequence has a further sub-subsequence 
  in which (*) holds, and -- up to a change in notation -- it is enough to show that the original sequence has some subsequence 
   in which (*) holds.
   
   Consider the set of all possible subsequential weak limits of sequences $T_{m_n}(\omega_n)$.  This is compact, so has 
at least one element $\mu$ which is maximal with respect to the ``stochastic order" partial order.  And a subsequence 
   $T_{m_n}(\omega^*_n)$ converging to $\mu$ clearly has property (*), because for any $\omega_n \le \omega^*_n$ 
   we have $T_{m_n}(\omega_n) \ge T_{m_n}(\omega^*_n)$ and so  by maximality $T_{m_n}(\omega_n)$ must also converge to $\mu$, as must $T_{m_n}(\sfrac{1}{2} \omega_n)$.

 \section{Discussion}
 
 \subsection{Regarding assumption (\ref{hypo})}
 Consider the ``line" graphs with
$\bV_n = \{1,2,\ldots,n\}$ and $w_{i,i+1} = 2^{-i}$.  
Here it is not possible to rescale time so that  assumption (\ref{hypo}) holds.
Heuristically, failure of assumption (\ref{hypo}) relates to this kind of exponential slowdown of edge-weights 
at the time of formation of the incipient giant component.

 \subsection{An example}
 \label{sec:example}
 Let us replace the assumption (\ref{hypo}) by the assumption
  \begin{equation}
   \lim_{t \to 0}  \limsup_n  \Ex C_n(t)/|\bV_n| = 0; 
   \quad
 \lim_{t \to \infty}  \liminf_n  \Ex C_n(t)/|\bV_n| = 1.
 \label{hypo-2}
 \end{equation}
 This is essentially saying, via a compactness argument, that the processes 
 $(C_n(t)/|\bV_n|, 0 \le t < \infty)$ converge to a limit process
  $(\widetilde{C}_\infty(t), 0 \le t < \infty)$ for which $\lim_{t \to 0} \widetilde{C}_\infty(t) = 0$ and $\lim_{t \to \infty} \widetilde{C}_\infty(t) = 1$.
 
Consider, for fixed $0<s<1$, the random times
 \begin{equation}
 T^*_n(s)  := \inf \{t: C_n(t) \ge  s  |\bV_n| \} .
 \label{T*n}
 \end{equation}
 Assumption (\ref{hypo-2}) implies these times 
 are $O(1)$ as $n \to \infty$, but is  {\em not} sufficient to show the
  ``weak concentration" property
\begin{equation}
\var \ T^*_n(s) \to 0 \mbox{ as } n \to \infty ,
\label{w-c}
\end{equation}
as the following example shows.

For the complete graph on $m$ vertices with edge-weights $1/m$, our bond percolation process is essentially just the 
 {E}rd{\H o}s-{R}\'enyi  process $\GG(m,t/m)$, for which the limit of the process
$(C_m(t)/m, 0 \le t < \infty)$ is a certain continuous function $\theta(\cdot)$ with $\theta(t) = 0$ for $t \le 1$ and
$\theta(t) > 0$ for $t > 1$ (explicitly,  $\theta(t)$ is the solution of $1 - \theta = \exp(-t \theta)$ -- see e.g. 
\cite{alon-spencer} sec. 10.4).
 Now take two copies of that complete graph on $m$ vertices with edge-weights $1/m$, and add a single edge $e^*$ between them with weight $m$.  This gives a graph on $n = 2m$ vertices.  
 It is easy to see that the limit of the process 
  $(C_n(t)/n, 0 \le t < \infty)$ is the random process
 \begin{eqnarray}
 \widetilde{C}_\infty(t) &=& \theta(t)/2, \ 0 \le t < \zeta \label{Cz} \\
 &=& \theta(t), \ \zeta \le t < \infty \nonumber
 \end{eqnarray}
 where 
$ \Pr(\zeta \le t) = \theta^2(t)$; 
here $\zeta$ represents the first time at which the giant component in each half contains the end-vertex of $e^*$.  
So the  ``weak concentration" property (\ref{w-c}) does not hold.

Note that we can modify this construction to make the graph {\em vertex-transitive}, that is there is a graph automorphism that
 maps any vertex to any other vertex.  Instead of a single edge between the two original copies, we assign weight $1/m^2$ to {\em every} 
 such edge.  
 Now the limit is again of form (\ref{Cz})  where now 
   $ \Pr(\zeta \le t) = 1 - \exp( -t \theta^2(t))$.

 \subsection{Analogies with bond percolation on infinite graphs}
 Rigorous mathematical treatment of bond percolation has focussed on infinite graphs, with ``general theory" developed under the assumption of {\em transitivity}, that is spatial symmetry. 
As the 2006 survey \cite{haggstrom} says,
\begin{quote}
\ldots infinite graphs,  where  the  issue  of  uniqueness  of  the  giant
component translates naturally into the question of whether there is a unique
infinite cluster.  This  has  the  advantage  of  always  having  a  clear-cut  yes/no 
answer, in contrast to the finite setting where it is not always totally obvious
what one really should mean by a giant component.
\end{quote}
In our setting of a {\em sequence} of finite edge-weighted graphs, one can readily formalize the idea 
 of giant components being unique as the property
 \begin{equation}
 \sup_t C^{[2]}_n(t)/|\bV_n| \to_p 0 \mbox{ as } n \to \infty 
 \label{UGC}
 \end{equation}
 where $C^{[2]}_n(t)$ is the size of the {\em second}-largest component at time $t$.
 Note that this can be restated in terms of the jumps of $(C_n(\cdot))$, as
  \begin{equation}
 \sup_{0 \le t < \infty}  | C_n(t) - C_n(t -) | / |\bV_n| \to_p 0 \mbox{ as } n \to \infty .
 \label{conj-assert}
 \end{equation}

\subsection{Two conjectures}
\label{sec:conjs}
Under the background assumption (\ref{hypo-2}), what further assumptions might be sufficient to imply either 
the ``unique giant component"  property (\ref{UGC}) or the ``weak concentration" property (\ref{w-c})? 
The example at the end of section \ref{sec:example} shows that {\em vertex-transitive} is not sufficient for either property.
If instead we assume {\em edge-transitive} then all edge-weights are equal and so we are in the more familiar 
setting of bond percolation on an {\em unweighted} graph with symmetry.   
Here it seems likely that known methods used in the infinite setting will be relevant.
However this requires care:  the infinite $r$-regular tree  does not have the ``unique giant component"  property but  typical realizations of
random $r$-regular graphs, as $n \to \infty$, do have this property \cite{asaf}, even though their local weak limit is the infinite $r$-regular tree.
Here is a bold conjecture.
  \begin{Conjecture}
  \label{C1}
  Consider a sequence of edge-transitive graphs with $|\bV_n| \to \infty$. 
  Then we can always rescale the edge-weight so that  (\ref{hypo-2}) holds.  After such rescaling, 
  the ``weak concentration" and the ``unique giant component"  properties hold.
  \end{Conjecture}
A second bold conjecture is that, without any assumption of symmetry, one of these properties implies the other.
  \begin{Conjecture}
Under   assumption (\ref{hypo-2}),   the ``unique giant component"  property (\ref{UGC}) implies the 
``weak concentration" property (\ref{w-c}).
\end{Conjecture}
Essentially, the conjecture is saying that the limit process $C_\infty(\cdot)$ indicated at the start of section \ref{sec:example} 
 might be deterministic and continuous (as in the classical settings) or might be random and discontinuous (as in the examples in section \ref{sec:example}), but cannot be 
 random and continuous.
(Readers aware of the famous
open problems involving continuity of the percolation function on infinite graphs should note that the  ``weak concentration" property 
relates to the  {\em inverse} of that function).

Separate from the literature on scaling windows mentioned in section \ref{sec:intro}, there is a line of work including
\cite{alon-benjamini,MR2917769} on bond percolation for unweighted finite graphs under isoperimetry assumptions, that is for expanders, which includes 
results on uniqueness of giant component.  Conjecture 1.1 of \cite{alon-benjamini}, not involving isoperimetry assumptions, is somewhat similar to our Conjecture \ref{C1}.

\subsection{Are there analogous results for first passage percolation?}

 Our starting structure was a  finite connected graph $(\bV,\bE)$ with edge-weights $\bw = (w_e)$ and with
independent Exponential(rate $w_e$) random variables $\xi_e$ associated with the edges.  
This structure can alternatively be used to construct first-passage times 
$X(v,v^\prime)$, defined as the minimum of $\sum_{e \in \pi} \xi_e$ over all paths $\pi$ from $v$ to $v^\prime$.
Regarding this as a model for spread of infection from an initial site $v$, 
the set of infected sites at time $t$ is
\[ \SS(v,t) := \{v^\prime: X(v,v^\prime) \le t \} . \]
Write
\[ \Delta :=  \max_{v,v^\prime} \Ex X(v,v^\prime) .\]
Given a sequence of such graphs with $|\bV_n| \to \infty$ and
a sequence $v_n \in \bV_n$,
consider the times 
 \begin{equation}
  T_n(v_n,s) :=  \inf \{t: \ | \SS_n(v_n,t) |  \ge s  |\bV_n| \} ; \quad 0 < s < 1 
    \label{def-Tnvs}
  \end{equation}
 and the ``incipient pandemic" time
  \begin{equation}
  T_n(v_n) := \inf \{t: \ | \SS_n(v_n,t) |  \ge |\bV_n|/\omega_n\} .
  \label{def-Tnv}
  \end{equation}
There is a simple analysis of such times,  provided we impose another assumption.
  Setting
  \[ w^* := \min \{w_e \ : \ w_e > 0\} \]
  we have
    \begin{Lemma}
 \label{Pvar}
 $ \var \ X(v,v^\prime) \le \Ex X(v,v^\prime)/ w^* $.
  \end{Lemma}
 Bounds of this type are classical on $\Ints^d$  \cite{kesten-speed} and are at least folklore in more general settings: an
 explicit statement and martingale proof in our setting is given in \cite{me-FPPG}. 
But inspecting the proof of Proposition 7 in  \cite{me-FPPG} shows that the same bound 
\begin{equation}
  \var \ T_n(v_n,s) \le \Ex T_n(v_n,s) / w^*_n 
  \label{varT}
  \end{equation}
holds for $T_n(v_n,s)$ at (\ref{def-Tnvs}), for arbitrary $s$. 

Now assume that in a sequence of graphs
\begin{equation}
w^*_n \Delta_n \to \infty .
\label{wnd}
\end{equation}
Because $\Ex T_n(v_n,s) = O(\Delta_n)$ for fixed $s$, (\ref{varT}) and (\ref{wnd})  imply the ``weak concentration" property
\[ \frac{T_n(v_n,s)}{\Delta_n} -  \frac{\Ex T_n(v_n,s)}{\Delta_n} \to_p 0 .\]
Compactness arguments as in section \ref{sec:compact} then lead to
a conclusion analogous to Proposition \ref{P1} for the ``incipient pandemic" time:
\begin{Corollary}
\label{C-pan}
Under assumption (\ref{wnd}),
there exists a deterministic sequence $\tau_n(v_n) \in [0,1]$ such that, for every 
 sequence  $\omega_n \uparrow \infty$ sufficiently slowly, the random times
$T_n(v_n)$ at (\ref{def-Tnv}) satisfy
 \[ T_n(v_n)/\Delta_n - \tau_n(v_n) \to_p 0 . \]
 \end{Corollary}
 However, this result in the   
 first-passage percolation setting differs in  two respects from the Proposition \ref{P1} result in the bond percolation setting. 
 To make Corollary \ref{C-pan} interesting we want the sequence $\tau_n(v_n)$ to be bounded away from zero, which is tantamount to the assumption (analogous to the first part of  (\ref{hypo})) that
  for some $0 < t_1  < \infty$ 
 \begin{equation}
 \lim_n  \Ex  | \SS_n(v_n,t_1 \Delta) | /|\bV_n| = 0 .
 \label{hypo-3}
 \end{equation}
 But in classical settings such as nearest-neighbor first-passage percolation on $\Ints^d_m$ \cite{kesten-FPP,50FPP}  this does not hold, 
 because by the shape theorem the scaling limit of $ | \SS_n(v_n,t \Delta) | /|\bV_n| $ is a deterministic function 
 $\phi(t)$ with $\phi(t) > 0$ for $t > 0$.
 The context where we do expect (\ref{hypo-3}) to hold is where the epidemic starts with 
 faster than polynomial growth, for instance on expander graphs or familiar models of random graphs \cite{shankar-1,shankar-2}. 
  Second, while assumption (\ref{wnd}) is stronger than necessary,
  we do need some assumption to prevent $T_n(v_n)$ have variability due to the influence  of a single edge-traversal time $\xi_e$ 
  associated with a very small weight $w_e$.
  For weak concentration of point-to-point percolation times $X(v,v^\prime)$,
  precise conditions in terms of such influence are given in \cite{me-FPPG}.
  It seems plausible that there are analogous precise conditions for weak concentration of the ``incipient pandemic" time $T_n(v_n)$,
  but we have  not studied this issue.

\paragraph{Acknowledgements.} 
I thank Asaf Nachmias for helpful comments.


 \end{document}